\documentclass{amsart}
\input xy
\xyoption{all}

\usepackage{amsfonts}                   

\begin{document}

\newcommand{\C}{{\mathbb C}}
\newcommand{\cee}{{\mathcal C}}                        
\newcommand{\cob}{{\bf Cob}}
\newcommand{\dee}{{\mathcal D}}
\newcommand{\diff}{{\rm Diff}}
\newcommand{\HH}{{\mathbb H}}
\newcommand{\met}{{\rm Met}}
\newcommand{\pt}{{\rm pt}}
\newcommand{\rel}{{\; {\rm rel} \;}}
\newcommand{\Q}{{\mathbb Q}}
\newcommand{\R}{{\mathbb R}}
\newcommand{\SO}{{\rm SO}}
\newcommand{\SU}{{\rm SU}}
\newcommand{\Spin}{{\rm Spin}}
\newcommand{\T}{{\mathbb T}}
\newcommand{\Wh}{{\rm Wh}}
\newcommand{\Z}{{\mathbb Z}}

\title {Topological gravity in Minkowski space}\bigskip
\author{Jack Morava}
\address{Department of Mathematics, Johns Hopkins University,
Baltimore, Maryland 21218}
\email{jack@math.jhu.edu}
\thanks{The author was supported in part by the NSF}
\date {20 May 2006}

\begin{abstract}
The two-category with three-manifolds as objects, $h$-cobordisms as
morphisms, and diffeomorphisms of these as two-morphisms, is extremely 
rich; from the point of view of classical physics it defines a nontrivial 
topological model for general relativity. \bigskip

\noindent
A striking amount of work on pseudoisotopy theory [Hatcher,
Waldhausen, Cohen-Carlsson-Goodwillie-Hsiang-Madsen \dots] can be 
formulated as a TQFT in this framework. The resulting theory is
far from trivial even in the case of Minkowski space, when the 
relevant three-manifold is the standard sphere.
\end{abstract}

\maketitle

\noindent
Topological gravity [18] extends Graeme Segal's ideas about conformal 
field theory to higher dimensions. It seems to be very interesting,
even in {\bf extremely} restricted geometric contexts: \bigskip

\begin{center}
{\bf \S 1 basic definitions}
\end{center} \bigskip

\noindent
{\bf 1.1} A {\bf cobordism} $W: V_0 \to V_1$ between $d$-manifolds is a
$D = d+1$-dimensional manifold $W$ together with a distinguished
diffeomorphism
\[
\partial W \cong V_0^{op} \coprod V_1 \;;
\]
a diffeomorphism $\Phi : W \to W'$ of cobordisms will be assumed
consistent with this boundary data. \bigskip

\noindent
$\cob(V_0,V_1)$ is the category whose objects are such
cobordisms, and whose morphisms are such diffeomorphisms. Gluing
along the boundary defines a composition {\bf functor}
\[
\# : \cob(V',V) \times \cob(V,V'') \to \cob(V,V'') \;.
\]
The two-category with manifolds as objects and the categories $\cob$ as
morphisms is symmetric {\bf monoidal} under disjoint union. \bigskip

\noindent
The categories $\cob$ are topological {\bf groupoids} (all
morphisms are invertible), with classifying spaces
\[
|\cob(V_0,V_1)| = \coprod_{[W:V_0 \to V_1]} B \diff (W \rel \partial) \;.
\]
The {\bf topological gravity} category has these objects as hom-spaces:
it is a (symmetric monoidal) topological category. \bigskip

\noindent
{\bf 1.2} A {\bf theory of topological gravity} is a representation of such a category
in some simpler monoidal category, e.g. Hilbert spaces, or spectra. \bigskip

\noindent
The homotopy-to-geometric quotient map
\[
B\diff = \met  \times_{\diff} E \diff \to \met \times_{\diff}{\rm pt} = 
\met/\diff
\]
defines a functor from the topological gravity category to a category with
the spaces $\met/\diff$ as morphism objects; these are the spaces of states in
general relativity (and 
\[
g \mapsto \int R(g) \; d{\rm vol}_g
\]
is a kind of Morse function upon them). \bigskip

\noindent
In Segal's conformal field theory, the corresponding objects are moduli 
spaces of (complex structures on) Riemann surfaces. Indeed if $W = 
\Sigma$ is a Riemann surface of genus $> 1$, its group of diffeomorphisms 
is homotopically discrete: the map
\[
\diff(\Sigma) \to \pi_0 \; \diff (\Sigma) 
\]
is a homotopy equivalence. The mapping class group acts with finite
isotropy on Teichm\"uller space, so when $d=1$ the homotopy-to-geometric 
quotient is close to a rational homology equivalence. \bigskip

\begin{center}
{\bf \S 2  examples}
\end{center} \bigskip

\noindent
{\bf 2.1} In recent work Galatius, Madsen, Tillmann and Weiss have identified
the classifying space of the cobordism category of oriented $d$-manifolds in terms
of a twisted desuspension $MT\SO(D)$ of the classifying space of the special orthogonal 
group. Their techniques extend more generally, to cobordism categories of manifolds
with extra structure on their tangent bundle. \bigskip

\noindent
Three-manifolds under Spin cobordism have very interesting connections with 
the theory of even unimodular lattices [8,16], and the methods of [6] identify the
classifying spectrum of this category with the desuspension of $B\Spin(4)$ by the
vector bundle associated to the standard four-dimensional representation of the 
spin group. Because of well-known coincidences in low-dimensional geometry, 
$\Spin(4) \cong \SU(2) \times \SU(2)$, so we can identify its classifying space 
with the product of two copies of infinite-dimensional quaternionic
projective space, and the vector bundle defined by the standard representation 
with the tensor product (over $\HH$) of the resulting two canonical quaternionic
line bundles $L_\pm$ ; thus
\[
MT\Spin(4) \sim (\HH P_\infty \times \HH P_\infty)_+^{-L_+^{op} \otimes_{\HH} L_-} \;.
\] 
The generators of $\pi_0 \Omega^\infty MT\Spin(4) \cong \Z^2$ can be identified with 
the signature and Euler characteristic, or alternately with the number of hyperbolic
and $E_8$ factors in the middle-dimensional intersection form [2] of a spin cobordism. 
\bigskip

\noindent
{\bf 2.2} There are other extremely interesting variant constructions in dimension four: 
contact three-manifolds under $\Spin^c$ cobordism define a natural context for 
Seiberg-Witten theory, while Lorentz cobordism [20] incorporates an arrow of time; 
but this note is concerned with {\bf 3-manifolds up to $h$-cobordism:}
\bigskip

\noindent
Recall that $W : V_0 \to V_1$ is an $h$-cobordism if the two inclusions 
\[
V_0 \subset W, \; V_1 \subset W
\]
are homotopy equivalences [17]. \bigskip

\noindent
The {\bf trivial} $h$-cobordism $W = V\times I$, where $I$ is an interval, 
is an interesting example. In dimensions $\geq 5$, the $s$-cobordism 
theorem classifies $h$-cobordisms by elements of the Whitehead group
\[
{\rm Im} \; [\pm \pi_(V) \to K_1(\Z[\pi_1(V)])] := {\rm Wh}(\pi_1(V)) \;,
\]
and there are invariants for {\bf parametrized} $h$-cobordisms taking
values in higher homotopy groups of certain pseudoisotopy spaces, which 
have been studied by Hatcher, Waldhausen, Igusa, \dots \bigskip

\noindent
This category has a monoidal structure, but it is relatively trivial, so that 
it is natural to assume that the manifolds $V$ are {\bf connected}. \bigskip

\noindent
{\bf 2.3} Here I will be concerned mostly with the case $V = S^3$: by Minkowski 
space I really mean the universal cover $S^3 \times \R$ of Penrose's (and 
others') conformal compactification $S^3 \times_{\pm 1} S^1$ of Minkowski 
space; this contains, in particular, a copy of Einstein's static universe [11]. 
Its time-like intervals define trivial $h$-cobordisms of $S^3$. \bigskip

\noindent
Note that there are lots of wild $S^3 \times \R$'s: remove a point from a 
fake $\R^4$. It would be very interesting to construct a semigroup of such 
things, under some kind of boundary gluing, as Segal did with topological 
annuli; current work of Gompf [10 \S 7, cf. also [3]] seems close to this. 
It is not clear at the moment if nontrivial smooth $h$-cobordisms of the 
three-sphere exist; the question is closely connected to the smooth 
four-dimensional Poincar\'e conjecture. \bigskip

\begin{center}
{\bf \S 3 double categories}
\end{center} \bigskip

\noindent
{\bf 3.1} Boundary value problems involve the interplay between diffeomorphisms of
a manifold and diffeomorphisms of its boundary. Tillmann [21] suggests that 
{\bf double} categories provide a natural framework for such questions. In
this context, the primary objects are certain rectangular diagrams 
\[
\xymatrix{
W : \ar@{=>}[d]_\Phi & V_0 \ar[d]^{\phi_0} \ar[r] & V_1 \ar[d]^{\phi_1}
\\ W': & V'_0 \ar[r] & V'_1 \;.}
\]
with cobordisms displayed horizontally, and diffeomorphisms (which 
preserve some boundary framing) presented vertically; these can be 
patched together in either direction. More recently, Getzler [7] has 
used manifolds, together with suitable (eg separating) codimension
one submanifolds, to define morphisms in such contexts; this seems
particularly suited to the {\it millefeuille} examples of Gompf, which 
(if I understand correctly) can be regarded as smooth $h$-cobordisms 
between topological, but not necessarily smooth, three-spheres. 
\bigskip

\noindent
{\bf 3.2} In any case, the double category $\dee$ of {\bf trivial} $h$-cobordisms 
between {\bf ordinary} three-spheres is already extremely interesting. 
I don't know  how to associate a topological category to a double category 
in general, but in this case pseudoisotopy theory defines an equivalence 
with the two-category
\[
\coprod [\{V\}/\cee(V)]
\]
having manifolds $V$ as its objects, and Cerf's group $\cee(V)$ [13 \S 6.2] of 
pseudoisotopies (regarded as a category with one object) as its category of 
automorphisms: \bigskip

\noindent
These pseudoisotopies are diffeomorphisms of the cylinders $V \times I$,
equal to the identity map on $V \times 0$. There is a fibration
\[
\diff( V \times I \rel \partial) \to \cee(V) \to \diff(V)
\]
of groups, and {\bf concordance}
\[
\Phi, \Psi \mapsto \Phi \; \# \; (\phi_1 \times 1_I) \circ \Psi 
\]
of pseudoisotopies defines a homomorphism 
\[
\cee(V) \times \cee(V) \to \cee(V) \;.
\]
The classifying space $B\cee(V)$ is thus a monoid, and the topological 
category associated to this rectification of $\dee$ defines an {\it ad hoc} 
topologification (with one object for each $V$, and the topological monoid 
$B\cee(V)$ for its space of endomorphisms). The classifying $B^2 \cee(V)$ space of 
{\bf that} topological category is the totalization of the bisimplicial space 
defined by the category of trivial $h$-cobordisms of $V$. \bigskip

\noindent                    
{\bf 3.3} There is a natural stabilization map from $B^2 \cee(V)$ to Waldhausen's 
ring spectrum $A(V)$. In the language of TQFT's, this defines a functor 
from the gravity category of trivial $h$-cobordisms of $V$ to the category 
with $\{V\}$ as its object, and the group ring $S^0[\Omega \Wh^d(V)]$ as its
endomorphism object. [The map from $\Omega B^2\cee(V)$ to $\Omega A(V)$ factors
through the space ${\rm HCobord}^d(V) \sim \Omega \Wh^d(V)$ of stabilized 
$h$-cobordisms of $V$ [22].] This reveals Whitehead torsion (regarded as
an element of $\Z[\Wh]$) as perhaps the primordial example of a TQFT! \bigskip

\noindent
Note that Cerf's maps define a fibration
\[
B\diff(V) \to B^2\diff(V \times I \rel \partial) \to B^2 \cee(V)
\]
which looks like a presentation of this {\it ad hoc} classifying space for a 
double category as a fibration 
\[
|{\rm Vertical}| \to |{\rm Horizontal}| \to |{\rm Double}|
\]
built from classifying spaces for its component (vertical and horizontal) 
morphisms; but I don't know enough about double categories to guess if 
this might be an instance of something more general. \bigskip

\begin{center}
{\bf \S 4 about} $A(S^n)$ 
\end{center} \bigskip

\noindent
{\bf 4.1} Through the efforts of many researchers, a great deal is known about 
the algebraic $K$-theory of spaces; in particular, if $X$ is simply
connected (and of finite type) its $A$-theory can be calculated (at
least $p$-locally [4 \S 1.3]) from the topological cyclic homology
[14 \S 7.3.14] of $S^0[\Omega X]$. \bigskip

\noindent
Since this pretends to be a paper about physics, however, I will be 
content with some remarks about $A_*(X) \otimes \Q$, which is accessible
in more elementary terms. [I want to record here my thanks to Bruce Williams 
and Bjorn Dundas for walking me through a great deal of literature in this 
field, without suggesting that they bear any responsiblity for the excesses 
of this paper.] \bigskip

\noindent
{\bf 4.2} In particular, old results [12] of Hsiang and Staffeldt imply that (when
$n > 1$) the rationalization of $A(S^n)$ splits as a copy of $A({\rm pt})_\otimes \Q
\cong K^{\rm alg}(\Z) \otimes \Q$ and the suspension of what is essentially the
(reduced) topological cyclic homology of $S^n$, which can be computed effectively
as the abelianization of $\tilde{H}_*(\Omega S^n,\Q)$ regarded as a graded Lie 
algebra; hence
\[
\pi_*(\Omega \Wh^d(S^n)) \otimes \Q \cong K^{\rm alg}_{*+1}(\Z) \otimes \Q \oplus 
\tilde{H}_*(\Omega S^n,\Q)_{\rm ab} \;.
\]
The Whitehead product structure on a wedge of spheres is rationally free, so the
graded Lie algebra structure has nontrivial commutators only when $n$ is even. 
When $n=2m+1$ is odd, the rational homology is polynomial on a single generator
$x_{2m}$; it follows that 
\[
A_{*+1}(S^3) \otimes \Q = \Q \langle \zeta_k, x_2^l \rangle
\]
is spanned as a rational vector space by elements $x_2^l$ of degree $2l$ and 
elements $\zeta_k$ of degree $4k$ corresponding to the odd zeta-values $\zeta(2k+1)$
which appear as regulators in Borel's calculations of $K^{\rm alg}_{4k+1}(\Z) \otimes \Q$.
\bigskip

\noindent
This can be made more precise; when $X$ is simply-connected then a reduced version 
$\Omega \widetilde{\Wh}(X)$ of loops on the Whitehead space is closely connected to 
a similarly reduced version $Q(\widetilde{LX}_{h\T})$ of (the infinite loopspace defined
by) the suspension spectrum of the homotopy quotient (by its natural circle action) of 
the free loopspace of $X$. \bigskip
 
\noindent
{\bf 4.3} The construction $Q = \Omega^\infty \Sigma^\infty$ sends a space to the infinite 
loopspace representing its suspension spectrum: this sends the rational homology 
of a space to its symmetric algebra. The cohomological invariants defined by the space of
trivial $h$-cobordisms of the three-sphere thus resemble the `big' phase spaces [9] studied 
in quantum cohomology: for example, the stable rational homology of the Riemann 
moduli space is essentially with the symmetric algebra on the homology of $\C P^\infty$,
and is thus a polynomial algebra with one generator of each even degree. \bigskip

\noindent
The rational cohomology of the infinite loopspace $\Omega^{\infty +1} A(S^3)$ seems 
similar in many ways: it is again a polynomial algebra, now with one set of generators
indexed by even integers, the other by integers $\equiv 0$ modulo four. Physicists see 
these symmetric algebras as Fock representations associated to certain polarized symplectic 
vector spaces. In our context this seems to be related to an `almost' splitting
\[
HC_{\rm per} \sim HC \oplus {\rm Hom}_{\Q[u]}(HC,\Q[u]) \;,
\]
of periodic cyclic homology [5] These representations have symmetries closely related 
to the Virasoro algebra, which lead [19] to interesting integrable systems. \bigskip

\noindent
This connection between 4D topological gravity and the equivariant free loopspace
of the three-sphere resembles in many ways a purely mathematical instance of a 
phenomenon physicists [1] call `holography', in which one physical model on the 
interior of a manifold is described by some other model on its boundary. Rather
than proceed any further with speculations along these lines, I'd like to close 
by raising a mathematical question: \bigskip

\noindent
A trivial $h$-cobordism between three-spheres is an example of a four-dimensional
spin cobordism; this defines a monoidal functor, and hence a morphism
\[
\Sigma^{-1} A(S^3) \to MT\Spin(4) \sim (\HH P_\infty \times \HH P_\infty)_+^{-L_+^{op} \otimes_{\HH} L_-} 
\]
of spectra. Could it possibly be nontrivial? \bigskip \bigskip

\bibliographystyle{amsplain}

\end{document}